\documentclass[preprint,12pt]{elsarticle}

\usepackage{graphicx}
\usepackage{subcaption}
\usepackage[labelformat=parens,labelsep=quad,skip=3pt]{caption}
\usepackage{graphicx}

\usepackage{amssymb}

\usepackage{lineno}

\usepackage{mathtools}

\journal{Integral Transforms and Special Functions}

\begin{document}

\begin{frontmatter}

\title{Finite Transforms with Applications to  Bessel Differential Equations of Order  Higher than Two}

\author{Gabriel L\'opez Garza}

\address{Mathematics Department, Universidad Aut\'onoma Metropolitana\\
Ciudad de M\'exico, M\'exico}

\begin{abstract}
A finite transformation method is introduced. This method is equivalent to the $Z$ transform method to a certain extent but generalizes it. By applying the presented method to the Bessel functions, it is possible to solve related ordinary differential equations of order higher than two with given initial conditions.
\end{abstract}


\end{frontmatter}


\section{Introduction}

The best-known transformation defined by sequences, among  the so-called ``Finite Transforms'', is probably the $Z$ transform (for an extensive study see \cite{V}).
Given a sequence $\{f_n\}$, the Z transform is defined by $F(z)=\sum_{n=0}^\infty\frac{f_n}{z^ n}$.
The main difference between the $Z$ transform with the transforms introduced in this paper is that whereas a sequence is transformable if the series in the definition converges for at least one value of $z$,  for the transforms that we study, convergence is not required. In section \ref{zeta} it is shown how the zeta transform is a special case of the transforms studied in this paper.
 Other well known examples of finite transforms are the Sturm-Liouville Transforms (see for instance \cite{CE}). 
 Basically, in order to solve some differential equations (ordinary or even partial differential equations), 
a sequence can be associated to certain suitable functions.
As it is known, the eigenfunctions $\{\phi_n(t)\}$ of self-adjoint operators satisfying certain boundary conditions (i. e., the solutions of a given Sturm-Liouville problem), provide a generalized Fourier expansion for a function $f$ (which satisfies certain conditions depending on the problem under consideration), given by
$$ f(t)=\sum_{n=0}^\infty a_n\phi_n(t),\quad a_n=\langle f,\phi_n\rangle,$$
where $\langle\cdot , \cdot \rangle$ is the inner product defined for each Sturm-Liouville problem. So, in this example of finite transform, the sequence $\{a_n\}$, may be defined as the transform, say $\mathcal{T}$, of a function $f(t)$ by
$$\mathcal{T}[f(t)]\stackrel{def}{=}\{a_n\},\,n\geq 0.$$
The transform so defined is very useful in solving differential equations to some extent. Still, the method induced by this transform (for instance in \cite{CE}) does not provide a complete  Mikusi\'nski's operational method as defined in \cite{M}. This is so since the convolution product induced by the inner product has non-zero divisors of zero (for instance for any $n\neq m,$ $\langle \phi_n,\phi_m\rangle =0$ given that the eigenfunctions of Sturm-Liouville problems are orthogonal and, of course, are different from zero). This limitation (limitation from the heuristic standpoint) is a rich source of examples for convolution algebras \cite{HH}, nevertheless.
The main interest of this paper is trying to extend the finite transform method as much as possible, if not to complete operational calculus in the Mikusi\'nski sense, but to extend it to a method that allows solving many ordinary or partial differential equations involving Bessel and other operators in various practical problems.

Examples of finite transforms are,  the Legendre Transform \cite{Ch1}, and the Lagerre transform \cite{McC}. Still, there are many other forms of associating a sequence to a given function besides the Sturm-Liouville transform, for instance, the Neuman series \cite[Chapter XVI]{W}. In our approach, a transform is constructed through a differential operator which is in some sense an extension of the Maclaurin series for Bessel Functions as will be explained.

The transform method that is studied in this paper consists in associating to a certain set of suitable functions a sequence $\{a_n\}\in\mathbb{C}$. Since the set of sequences is not an algebraic field, the study is restricted to sequences that are invertible for the Cauchy product. The fact that not each sequence different from zero is invertible is the main cause of the impossibility of extending our study to a complete Mikusi\'nski's Operational Calculus. So our study shares more similarities with the traditional Laplace transform Method than with the Mikusi\' nski's Operational Calculus as that studied in \cite{BL}.
In fact, to certain differential equations with given initial conditions, a polynomial may be associated via the finite transformation, so that an equation for the transform may be solved in terms of partial fraction decomposition. The partial fraction decomposition may be associated with known transforms or the Cauchy product of known transforms. Finally, the inverse transform is obtained with which the problem is solved, as usual.

The article is divided as follows: in section \ref{mathset} the operations used in the set of transformations are described and the finite transform used is defined (subsection  \ref{DTnn}).
In section \ref{examplesn} some higher order Bessel equations are solved and equivalence to Laplace transform is shown but is worth noticing that the transforms are defined by differential and not with integral operators as in the case of classic Laplace transform. 
In the final section, the striking similarities that have long been noted \cite[p. 136] {D} between the functions ber and bei with the functions cos and sin, respectively, are fully explained within the context of the transform method studied in this article.

\section{Mathematical setting}
\label{mathset}
Given two sequences $\{a_n\},\{b_n\}$ we define the product (the so-called Cauchy product) by
\begin{eqnarray}
\label{Cp} \{a_n\}\{b_n\}&=&\left\{ \sum_{\tau=0}^n a_\tau b_{n-\tau}\right \}\\
\nonumber                                         &=& \{ a_0b_0,a_0b_1+a_1b_0,\dots\}.
\end{eqnarray}In \cite[p. 721]{Br} it is shown that the Cauchy product of two sequences is zero if and only if one of the sequences is the zero sequence $a_n=0,n\geq 0$ so that a quotient field of sequences can be constructed. Nevertheless not every sequence different from the zero sequence is invertible. In order to be invertible, the first term of a given sequence is required to be different from zero.
Actually, given a sequence $\{a_n\},a_n\in \mathbb{C}, a_0\neq 0$, the multiplicative inverse $\{ b_n\}$ with respect to the Cauchy product is easily calculated recursively,
since $\{a_n\}\{b_n\}=\{1,0,0,\dots\}$ implies
\begin{eqnarray*}
b_0&=&\frac{1}{a_0}\\
b_1&=&-\frac{a_1}{a_0^2}\\
b_2& = & \frac{a_1^2}{a_0^3}-\frac{a_2}{a_0^2}\\
      &\vdots &  \\
b_n & = & \frac{-1}{a_0}\left( a_1b_{n-1}+\cdots+a_n b_0\right).
\end{eqnarray*}
The fact that not every sequence different from zero has a multiplicative inverse restricts the construction of an operational calculus in the Mikusi\'nski's sense \cite{M}, but, as we will see, from the invertible sequences used in this paper, it is possible to build many transforms and operate with them as in many other transform methods used to solving differential equations.
After the last considerations the set of suitable sequences $\mathcal{A}=\{ \{a_n\}: a_0\neq 0\}$ is defined.

 With the Cauchy product, it is possible to construct many operators $T:\mathcal{A}\to \mathcal{A}$. An important example is 
the {\it right shift}, which is an operator defined by the Cauchy product with the sequence  $s=\{0,1,0,0,\dots\}$ if $\{a_n\},n\geq 0$ is any sequence with $a_n\in \mathbb{C}$, we have
$$s\{a_0,a_1,a_2,\dots\}=\{0,1,0,0,\dots\}\{a_0,a_1,a_2,\dots\}=\{0,a_0,a_1,a_2,\dots\}.$$
So the right shift operator $S:\mathcal{A}\to \mathcal{A}$ is defined by $S[\{a_n\}]\stackrel{def}{=}s\{a_n\}$. The notation $S[\{a_n\}]$ is not in use in most papers and calling the operator $s$ instead $S$ is the standard procedure (for instance in \cite{Br}), we will follow this practice throughout this paper to avoid confusion with an established practice.

It follows that $s^2=s s=\{0,0,1,0,0,\dots\}$ and, in general $$s^n=\{0,0,\dots,0,1,0,\dots\},$$
is the sequence with zeros everywhere except in the place $n+1$ where there is the number one. With the powers of $s$ the following notation is standard
\begin{eqnarray}
\label{sum}\{a_0,a_1,a_2,\dots\}\stackrel{def}{=}a_0s^0+a_1s+a_2s^2+\cdots,
\end{eqnarray}
where $s^0=\{1,0,0,\dots\}$ and subsequently we write simply $a_0s^0=a_0$, and any constant $c$ in our calculus represents the sequence $c\stackrel{def}{=}\{c,0,0,\dots\}$. We recall that the relation (\ref{sum}) is purely formal and does not involve the concept of convergence at all \cite[p. 723]{Br}.

 Notice that $s=\{0,1,0,\dots\}$ is not invertible with respect to the Cauchy product, but the sequence $\{1,0,0,\dots\}+s=\{1,1,0,0,\dots\}$, actually is, and the following  formulas are derived easily  \cite{Br}
\begin{eqnarray}
\label{tf14} \{1\} &=&\frac{1}{1-s},\\
\label{tf22} \{r^n\}& = & \frac{1}{1-rs},\\
\label{tf25} \left\{\cos \frac{\pi}{2}n\right\}&=&\{1,0,-1,0,\dots\}=\frac{1}{1+s^2},\\
\label{tf26} \left\{\sin \frac{\pi}{2}n\right\}&=&\{0,1, 0,-1,\dots\}=\frac{s}{1+s^2}.
\end{eqnarray}
That is, the sequence $1-s$ is the multiplicative inverse of the sequence $\{1\}=\{1,1,1,\dots\}$ for the Cauchy product or, more properly speaking, the sequence $1-s$ is a representative of the class of inverse sequences of $\{1\}$, and similar meaning has the symbol $1/G(s)$ for the identities from (\ref{tf22}) to (\ref{tf26}), as well.

The {\it left shift} $l$ of a sequence $\{y_n\}$ is defined by
\begin{eqnarray}
\nonumber l\{y_n\}=l\{y_0,y_1,y_2,\dots\}&\stackrel{def}{=}&\{y_1,y_2,y_3,\dots\}=\{y_{n+1}\}\\
\label{leftsm} l^m\{y_n\} &\stackrel{def}{=}& \{y_{n+m}\}.
\end{eqnarray} 

The following formula \cite[p. 724]{Br} will be relevant for solving differential equations related to Bessel and other operators in this paper:
\begin{eqnarray}
\label{fun} s^m(l^m \{y_n\})= s^m\{y_{n+m}\}=\{y_n\}-y_0-sy_1-\cdots- s^{m-1}y_{m-1}.
\end{eqnarray}
Notice that formula (\ref{fun}) is similar to the Laplace transform formula for the derivative of order $m$ of a given function.

\subsection{Finite Transform definition}
\label{DTnn}
Given a series $g(t)=a_0f_0(t)+a_1f_1(t)+a_2f_2(t)+\cdots$, where $a_n\in\mathbb{C}$ and $f_n(t)$ are given functions, we define the transform $\mathcal{T}[g(t)]=G(s)$ by the formula
\begin{eqnarray}
\label{ftdef} \mathcal{T}[g(t)]&\stackrel{def}{=}&\{a_n\}=\{a_0,a_1,a_2,\dots\}\\
 \label{ftdef2}                                                   &= &a_0+a_1s+a_2s^2+\cdots=G(s),
\end{eqnarray}
if and only if there exists a  operator $L$ such that $L[f_n(t)]=f_{n-1}(t), n\geq 1$. The operator $L$ corresponds to the left shift $s$ and it is called the {\it concrete realization} of the shift for the sequence $\{ f_n(t)\},n\geq 0$.

The $a_n$ will be given by different instances of differential operators applied to given functions as seen in the following examples.
Notice that (\ref{ftdef2}) is a purely formal expression in powers of $s$ and this notation is not germane with the question of convergence.

\section{Examples of applications}
\label{examplesn}

\begin{enumerate}
\item \label{ex1ex} {\it Transform induced by Bessel functions.} For Bessel functions of order $\nu\geq 0$ we consider the monomials of the form $f_{n,\nu}(t)=\frac{(t/2)^{2n+\nu}}{\Gamma (\nu+n+1)n!}$. If $L_\nu\stackrel{def}{=}\frac{1}{t}DtD-\frac{\nu^2}{t^2}$ is a differential operator where $D$ denotes the derivative with respect to $t$, then a direct calculation shows that
\begin{eqnarray}
\label{gg6} L_\nu f_{n,\nu}(t)&=&f_{n-1,\nu}(t),\text{ for }n\neq 0,\\
\label{gg7} L_\nu f_{0,\nu}(t)&=& 0, \text{ for }n=0,
\end{eqnarray}
in this way  $L_\nu $ corresponds to the left shift $s$ in this concrete realization.
In fact, if $f(t)$ is a function for which
$$f(t)=a_0f_{0,\nu}(t)+a_1f_{1,\nu}(t)+a_2f_{2,\nu}(t)+\cdots, a_n\in\mathbb{C},$$
then we define the Bessel Transform $\mathcal{T}_{B_\nu}$ by
$$\mathcal{T}_{B_\nu}[f(t)]=\{a_0,a_1,a_2,\dots\}$$
and, since
 $$L_\nu f(t)=a_1f_{0,\nu}(t)+a_2f_{1,\nu}(t)+\cdots,$$
by (\ref{gg6}) and (\ref{gg7}), we have
$$\mathcal{T}_{B_\nu}[L_\nu f(t)]=\{a_1,a_2,a_3\dots\}$$
so that applying $L_\nu$ to a function $f(t)$ corresponds to the left shift applied to $\mathcal{T}[f(t)]$, as we claimed.


Examples of transforms for fixed $\nu\geq 0$ are
\begin{eqnarray}
\label{bet1} \mathcal{T}_{B_\nu}[J_\nu (t)] &= &\{1,-1,1,-1,\dots\},\\
\label{bet2} \mathcal{T}_{B_\nu}[I_\nu (t)] & =& \{1,1,1,\dots\},\\
\label{bet3} \mathcal{T}_{B_\nu}[\text{Ber}_\nu (t)]& = &\left \{  \cos\frac{(3\nu+2n)\pi}{4} \right\},\\
\label{bet4} \mathcal{T}_{B_\nu}[\text{Bei}_\nu (t) ]& = & \left\{\sin \frac{(3\nu+2n)\pi}{4}\right\}.
\end{eqnarray}
where $J_\nu(t)$ are the well known Bessel functions of order $\nu\geq 0$ and $I_\nu(t)$ are the modified Bessel functions of the first kind. 
Particular interesting cases of (\ref{bet3}) and  (\ref{bet4}) occur when $\nu=0$ for which
\begin{eqnarray}
\label{bet5}  \mathcal{T}_{B_0}[\text{Ber} (t)]& = &\{1,0,-1,0,\dots \}\\
\label{bet6}  \mathcal{T}_{B_0}[\text{Bei} (t)]& = &\{0,-1,0,1,\dots \}.
\end{eqnarray}
Formulas (\ref{bet1}) to (\ref{bet6}) are well-known (see for instance \cite{S} p. 140 for formulas (\ref{bet3}) and (\ref{bet4})). For general functions, the Transform $\mathcal{T}_{B_\nu}$ can be calculated by iterating the operator $L_\nu$. Let  $L_\nu^m =L_\nu(L_\nu^{m-1}),m\geq 2$ be the order $m$ operator, then given a function $f(t)$ for which the limit
\begin{eqnarray}
\label{limB}\lim_{t\to 0}(L^mf(t))= a_m
\end{eqnarray}
does exist we will denote $\lim_{t\to 0}(L_\nu^m f(t))=L^m_\nu f(0)=a_m$, so for the class of functions for which that limit exists for any $m\in\mathbb{N}$ we define
$$\mathcal{T}_{B_\nu}[f(t)]=\{ f(0),L_\nu f(0),L_\nu^2f(0),\dots \}.$$
Note the similarity of the series $f(x)=f(0)+a_1f_{1,\nu}(t)+a_2 f_{2,\nu}(t)+\cdots$ with the Maclaurin series. The coincidence is not developed further, for our purposes, the existence of the transform $\mathcal{T}_{B_\nu}$ requires only the existence of the limit (\ref{limB}). There are many functions for which the transform can be computed by applying the $L_\nu$ operator repeatedly, for example, $\sqrt{\frac{2}{\pi t}} \sinh t$, $\sqrt{\frac {2}{\pi x}} \cosh t$, $\sqrt{\frac{2}{\pi t}} \left( \sin t +\frac{\cos t}{t} \right) $, and many others. The reader can verify directly by applying the operator that the transform $ \mathcal{T}_{B_\nu}$ exists, but a direct proof follows easily from the formulas for Bessel functions of order equal to one half and odd integer \cite[ pg. 138 and p. 140]{S}. For example, with the formula 24.58 in \cite{S} we have
$I_{1/2}(t)=\sqrt{\frac{2}{\pi t}}\sinh t$, and therefore

$$\mathcal{T}_{B_{1/2}}\left[\sqrt{\frac{2}{\pi t}}\sinh t \right]=\{1,1,1,\dots  \}.$$


\item\label{ex1bis} {\it Transforms method for Bessel operators.} After separation of variables of the Plum equation $\Delta^2u-\gamma \Delta u-\frac{4\gamma}{r^2}u=\Lambda u$, where $\Delta$ is the laplacian in polar coordinates \cite{Ev}, the resulting fourth order equation 
$ (ty'')''-((9t^{-1}+8\mu^{-1}t)y')'=\Lambda ty$ 
may be solved by using the transforms of example \ref{ex1ex}. 
In fact, last equation may be written in the form equivalent to \cite[eq. 45]{BL} 
\begin{eqnarray}
\label{eq1ex1bis}
\left[L^2_2-\frac{8}{\mu} L_2-\left( \lambda^2+\frac{8}{\mu}\right)\lambda^2\right]y(t)=0.
\end{eqnarray}
Applying $\mathcal{T}_{B_2}$ to equation (\ref{eq1ex1bis}) we obtain
\begin{eqnarray}
\label{teq1ex1}
\{ Y_{n+2} \}-\frac{8}{\mu} \{  Y_{n+1}\}-\left( \lambda^2+\frac{8}{\mu}\right)\lambda^2 \{  Y_n \}=0,
\end{eqnarray}
where $\{Y_n\}=\mathcal{T}_{B_2}[y(t)]$, $\mathcal{T}_{B_2}[L^2_2y(t)]$, and $\mathcal{T}_{B_2}[L_2y(t)]$  are obtained after using formula (\ref{leftsm}).
Multiplying (\ref{teq1ex1}) by $s^2$ and applying formula (\ref{fun}) we have
\begin{eqnarray*}
\{Y_n\}-Y_0-sY_1-s\frac{8}{\mu}(\{ Y_n\}-Y_0)-s^2\left( \lambda^2+\frac{8}{\mu}\right)\lambda^2\{ Y_n\}=0,
\end{eqnarray*}
and hence
\begin{eqnarray}
\label{pfB} \{ Y_n \}=\frac{Y_0+s\left(  Y_1-\frac{8}{\mu} Y_0\right)}{(1+\lambda ^2 s)\left(1-\left(\frac{8}{\mu}+\lambda^2\right)s\right)}
\end{eqnarray}
and after partial fraction decomposition in (\ref{pfB})
\begin{eqnarray}
\label{pfB1} \{ Y_n \}=\frac{1}{2(4+\lambda^2\mu)} \left(\frac{8Y_0-\mu Y_1+Y_0\lambda^2\mu}{1+\lambda^2 s}
                                   +\frac{Y_1\mu+Y_0\lambda ^2\mu}{1-(8/\mu+\lambda^2)s}\right).
\end{eqnarray}
By using  a general form of formulas (\ref{bet1}) and (\ref{bet2}),  i. e.
\begin{eqnarray}
\label{13g} \mathcal{T}_{B_\nu}[J_\nu(\sqrt{\lambda}t)]& = & \lambda^{\nu/2}\{\lambda^k\}=\frac{\lambda^{\nu/2}}{1-\lambda s}\\
\label{14g}  \mathcal{T}_{B_\nu}[I_\nu(\sqrt{\lambda}t)]& = & \lambda^{\nu/2}\{(-\lambda)^k\}=\frac{\lambda^{\nu/2}}{1+\lambda s}
\end{eqnarray}
we have
\begin{eqnarray*}
\label{bet12n}\mathcal{T}_{B_2}[J_2(\lambda t)]&=&\lambda^2\{1,-1,1,-1,\dots\}=\frac{1}{1-(-\lambda^2)s}\\
\label{bet22n}\mathcal{T}_{B_2}\left[I_2\left(\sqrt{\frac{8}{\mu}+\lambda^2}\, t\right)\right]&=&\left(\frac{8}{\mu}+\lambda^2 \right)\{1,1,1,1,\dots\}=\frac{1}{1-\left(\frac{8}{\mu}+\lambda^2 \right)s}. 
\end{eqnarray*}
With the last two equations, it is possible to take inverse transform in equation (\ref{pfB1}) to solve the initial value problem (\ref{eq1ex1bis}) with initial conditions
$y(0)=Y_0$, and $\lim_{t\to 0}L_2y(t)=Y_1$.
The reader may notice that with the identities $J_2(u)=-J_0(u)+\frac{2}{u}J_1(u)$ and $I_2(u)=I_0(u)-\frac{2}{u}I_1(u)$ we obtain the same solution as that given by formulas in \cite{Ev8} and in \cite{BL}.

\item\label{ex2ex} {\it Transform induced by Maclaurin series.} Consider the monomials of the form $f_n(t)=t^n/n!$, so that,  since $\frac{d}{dt}t^n/n!=t^{n-1}/(n-1)!$, the  left shift in this concrete realization is the derivative with respect to $t$. 
Given a 
function $f(t)$ 
for which $\lim_{t\to 0}f^{(n)}(t)=f^{(n)}(0),\,n\geq 0$ does exist, the discrete transform of $f$, say $\mathcal{T}_M[f(t)]=F(s)$ is, by definition, the sequence
\begin{eqnarray*}
\mathcal{T}_M[f(t)]&\stackrel{def}{=}&\{ f(0), f'(0), f^{(2)}(0),\dots, f^{(n)}(0),\dots\}\\
            &=&f(0)+f'(0)s+f^{(2)}(0)s^2+\cdots,
\end{eqnarray*}
and, formally,$$\mathcal{T}_M[f'(t)]=\{f'(0),f^{(2)}(0), f^{(3)}(0),\dots\}=f'(0)+f^{(2)}(0)s+f^{(3)}(0)s^2+\cdots.$$
In general if $\mathcal{T}_M[y(t)]=\{a_0,a_1,a_2,\dots\}=\{a_n\}, n\geq 0$, then $\mathcal{T}_M[y'(t)]=\{a_1,a_2,\dots\}=\{a_{n+1}\},n\geq 0$. Consequently $\mathcal{T}_M[y^{(m)}(t)]=\{a_m,a_{m+1},\dots\}=\{a_{n+m}\},n\geq 0$, so that the derivative of order $m$ corresponds to the right shift of $m$ places.

Observe that, for instance, the sequence $\{1,1,1,\dots\}$ corresponds to the exponential function, and by formula (\ref{tf14}), $\mathcal{T}_M[\,e^t\,]=\{1,1,1,\dots\}=\frac{1}{1-s}$.
 Also, it is easy to see, according to formula (\ref{tf25}), that $\mathcal{T}_M[\,\cos t\,]=\frac{1}{1+s^2}$. 
And, moreover, by formula (\ref{tf26}), $\mathcal{T}_M[\,\sin t\,]=\frac{s}{1+s^2}.$

With the formula  (\ref{fun}) now it is possible to solve non-homogeneous differential equations  with constant coefficients, as an example we solve 
\begin{eqnarray}
\label{ex21}y''-3y'+2y=e^{3t};\\
\label{ex22} y(0)=1,\;y'(0)=0.
\end{eqnarray}
 We have by (\ref{leftsm}) that if $\mathcal{T}_M[y(t)]=\{y_n\},n\geq 0$ then in this realization, necessarily the derivative satisfies
\begin{eqnarray}
\label{deriv} \mathcal{T}_M[y^{(m)}(t)] =l^m\{y_n\}=\{y_{n+m}\}.
\end{eqnarray}

Taking transforms in (\ref{ex21}) and setting $\mathcal{T}_M[y(t)]=\{y_n\}=Y(s)$, we have by formula (\ref{deriv})
\begin{eqnarray}
\label{ex23}\{ y_{n+2}\}-3\{y_{n+1}\}+2\{y_n\}=\frac{1}{1-3s}.
\end{eqnarray}
Multiplying (\ref{ex23}) by $s^2$, taking into account the initial conditions (\ref{ex22}) so that $\mathcal{T}_M[y(0)]=y_0=1$, $\mathcal{T}_M[y'(0)]=y_1=0$, hence by formula (\ref{fun}),  we obtain after simplification
\begin{eqnarray}
\label{ex24} \{y_n\}(1-3s+2s^2)&=&\frac{s^2}{1-3s}+1-3s\\
\label{ex244} \{y_n\} & = & \frac{1-6s+10s^2}{(1-3s)(s-1)(2s-1)}\\
\label{ex25} Y(s)& = & \frac{\frac{1}{2}}{1-3s}+\frac{\frac{5}{2}}{1-s}-\frac{2}{1-2s}\\
\label{ex26} y(t)&=& \mathcal{T}_M^{-1}[Y(s)]=\frac{1}{2}e^{3t}+\frac{5}{2}e^t-2e^{2t}.
\end{eqnarray}
Clearly (\ref{ex25}) is obtained after partial fraction decomposition of (\ref{ex244}), and  (\ref{ex26}) is obtained from (\ref{tf22}), since in this realization
$\mathcal{T}_M[e^{rt}]=\{r^n\}=\frac{1}{1-rs}$. Of course, as usual  in transform methods, if $\mathcal{T}_M[y(t)]=Y(s)$ then we define $y(t)\stackrel{def}{=}\mathcal{T}_M^{-1}[Y(s)]$.

\end{enumerate}


\section{The $Z$ transform case}
\label{zeta}

A correspondence between the $Z$ transform and the Maclaurin transform studied in this article is established now. For a function $G\in  \mathcal{C}^\infty$ a sequence $\{g_n\}$ may be defined by the formula 
\begin{eqnarray}
\label{defz} g_n\stackrel{def}{=}\lim_{x\to 0}\frac{1}{n!}\frac{d^n}{dx^n}G\left(\frac{1}{x}\right)
\end{eqnarray}
if and only if the limit does exist. With the sequence $\{g_n\}$, the zeta transform of the sequence is defined as
\begin{eqnarray}
\label{defzt} \mathcal{T}_{Z}[\{g_n\}]=F(z)=g_0+\frac{g_1}{z}+\frac{g_2}{z^2}+\cdots=\sum_{n=0}^\infty\frac{g_n}{z^n},
\end{eqnarray}
if the series (\ref{defzt}) converges for at least one $z\in\mathbb{C}$. Properties of the zeta transform are well known (see for instance \cite{V}). 
Notice that the correspondence
\begin{eqnarray}
\label{corr} z^{-n}\leftrightarrow s^n
\end{eqnarray}
establishes a one-to-one correspondence, hence an equivalence between the zeta transform and Maclaurin transform. 
The principal difference between these transforms is that the left shift $s$ is associated with a differential operator in all other transforms studied in this paper, (and in particular for the Maclaurin transform $\mathcal{T}_M$, in such a way that $s$ corresponds to the derivative), meanwhile by derivating $z^{-n}$ is not possible to find the ${g_n}$ from a given Laurent series $\sum_{n=0}^\infty\frac{g_n}{z^n}$  and evaluating at $z=0$, but it is possible by complex integration.

Another difference is, as already mentioned, the convergence of a given series is not determinant for the existence of the Maclaurin transforms. For instance, the sequence $\{1,1,1,\dots,\}$ has  $Z$ transform 
\begin{eqnarray}
\label{zt1}
\mathcal{T}_Z[\{1,1,\dots\}]=F(z)=1+\frac{1}{z}+\frac{1}{z^2}+\cdots=\frac{z}{z-1}
\end{eqnarray}
which exists only if $|z|>1,$ but 
\begin{eqnarray}
\label{ztm}
\mathcal{T}_M[\{1,1,\dots\}]=\frac{1}{1-s}
\end{eqnarray}
is well defined.
So, meanwhile  in (\ref{zt1}) the $Z$ transform $F(z)$ only exists for $|z|>1$,  formula (\ref{ztm}) indicates that the multiplicative inverse respect to the Cauchy product of the sequence
$\{1,1,1,\dots\, 1\dots\}$ is the sequence $\{1,-1,0,0,\dots,0,\dots\}=1-s$, as the reader may easily corroborate. 

\noindent {\bf Example [Equivalence with $Z$ transform].} As an illustration of the equivalence between the $Z$ transform and the Maclaurin transform we solve an  initial value problem:
\begin{eqnarray}
\label{e342} y_{k+1}-3y_k&=&4\\
   \nonumber   y_0&=&1.
\end{eqnarray}
 which is already solved by using $Z$ transforms in \cite[Example 3.42]{K}. By the properties of the $Z$ transform \cite[Chapter 3, section 3.7]{V} equation and initial condition in (\ref{e342}) is transformed in
\begin{eqnarray}
\label{e342t} Y(z)=\mathcal{T}_Z[\{y_n\}]&=&\frac{-2z}{z-1}+\frac{3z}{z-3}\\
\label{e342t2}                                                              & =&-2\sum_{n=0}^\infty\frac{1}{z^n}+3\sum_{n=0}^\infty\frac{3}{z^n}
\end{eqnarray}
where (\ref{e342t2}) is obtained by expanding in Laurent series (\ref{e342t}). In considering the correspondence (\ref{corr}) and applying it to (\ref{e342t2}) it is possible to find the Maclaurin transform of (\ref{e342}) given by



\begin{eqnarray*}
\label{342M} Y(s)=\mathcal{T}_M[\{y_n\}]&=&\frac{-2}{1-s}+\frac{3}{1-3s}\\
                 \{ y_k\}= \mathcal{T}^{-1}_M[Y(s)]&=&-2 \{1,1,\dots\}+3\{1,3,3^2,\dots\}
\end{eqnarray*}

So the difference problem (\ref{e342}) has the solution $y_k=-2+3^{k+1}$, as in \cite{K}.

Of course, problem (\ref{e342}) can be solved directly with the Maclaurin transform, as shown below, last argument was made only to emphasize the correspondence (\ref{corr}).
To solve (\ref{e342}) directly with the Maclaurin transform $Y(s)=\mathcal{T}_M[\{y_n\}]$ we have by the equation in problem (\ref{e342})
$$ \{ y_{k+1}\}-3\{y_k\}=4\{1,1,\dots\},$$
multiplying the last equation by $s$, applying formula (\ref{fun}), and taking $y_0=1$, which corresponds to the initial value in  problem (\ref{e342}), it is  obtain
\begin{eqnarray*}
s\{y_k\}-1-3s\{y_k\}&=&4s\{1,1,\dots\}\\
(1-3s)Y(s)&=&\frac{4s}{1-s}+1\\
Y(s)&=&\frac{4s}{(1-s)(1-3s)}+\frac{1}{1-3s},
\end{eqnarray*}
consequently, after partial fraction decomposition
$$Y(s)=\frac{-2}{1-s}+\frac{2}{1-3s}+\frac{1}{1-3s}=\frac{-2}{1-s}+\frac{3}{1-3s}.$$
So by taking inverse transform and using formula (\ref{tf22})
$$\{y_k\}=-2\{1,1,\dots\}+3\{3^k\}, k>0,$$
so, $y_k=-2+3^{k+1}, k>0$ and $y_0=1$ which coincides with the $Z$ transform solution given before.

\section{Conclusions}

 The striking similarity between ber and bei functions with cos and sin functions respectively has been noticed since time ago (see for instance \cite[p. 136]{D}) but, to the best of my knowledge, it has not been completely understood. In the approach of this article, Ber and Bei functions and cos and sin are in correspondence with the same transform (or sequence),  Ber and Bei for $L_\nu$ operator, and cos and sin for $d/dt$ operator respectively. Table (\ref{tablet}) shows the exact correspondence between these functions.

\begin{table}[!h]
\begin{center}
\begin{tabular}{c| l l l } 
 \hline
 Operator & function & Transform & associated sequence \\ [0.5ex] 
 \hline
 $L_\nu$ & Ber$(t)$ &  $\mathcal{T}_{B_\nu}[\text{Ber}_\nu\,t]=\frac{1}{1+s^2}$&   $\{1,0,-1,0,\dots\}$  \\ 
 $L_\nu$ &  Bei$(t)$  &  $\mathcal{T}_{B_\nu}[\text{Bei}_\nu\,t]=\frac{s}{1+s^2}$ & $\{0,-1,0,1,\dots\}$  \\
 $L_\nu$ &  $I_\nu (t)$  &  $\mathcal{T}_{B_\nu}[I_\nu t]=\frac{1}{1-s}$ &    $\{1,1,\dots\}$  \\
$D$  &  $\cos t$ & $\mathcal{T}_M[\cos t]=\frac{1}{1+s^2}$ &   $\{1,0,-1,0,\dots\}$ \\
 $D$  & $\sin t$ & $\mathcal{T}_M[\sin t]=\frac{s}{1+s^2}$ &  $\{0,-1,0,1,\dots\}$   \\
 $D$  & $e^t$ & $\mathcal{T}_M[e^t ]=\frac{1}{1-s}$ &    $\{1,1,\dots\}$   
\end{tabular}
\caption{\label{tablet} Correspondence between sin, cos, and Bessel functions according to their respective operator.}
\end{center}
\end{table}

So, for instance, Table (\ref{tablet}) shows that with sequence $\{1,1,\dots\}$ is in correspondence with two different functions, the Bessel function $I_\nu$ and the exponential function. 
So that they have the same transforms even when they were obtained with different transform methods. Of course, the operator associated with each function is $L_\nu$ and $D$, respectively.
So the coincidence of two different methods is fully explained within the context of the transform method studied in this article.
\bibliographystyle{model1-num-names}

\end{document}